
\magnification=\magstep1
\input amstex
\documentstyle{amsppt}
\hsize=6.5truein
\vsize=8.9truein
\NoBlackBoxes

\def\HH{\text{Hom}}
\def\H{\text {Hom}_R}
\def\Rad{\text{Rad}\thinspace}
\def\HO{\text{Hom}_R(\omega,\omega)}
\def\J{\bold j}
\def\CM{Cohen-Macaulay }

\def\HR{\hat R}
\def\HS{\hat S}
\def\Ho{\hat\omega}

\define\fraktur#1{\eufm#1}
\topmatter
\title INDECOMPOSABLE CANONICAL MODULES\\ AND CONNECTEDNESS \endtitle
\author Melvin Hochster and Craig Huneke \endauthor
\address Department of Mathematics, University of Michigan, Ann Arbor, MI 48109-1003, USA \endaddress
\address Department of Mathematics, Purdue University, West Lafayette, IN 47907 USA \endaddress
\endtopmatter
\document
\footnote""{1991 Mathematics Subject Classification. Primary 13E05, 13C99, 13D45, 13H99, 13J10.}
\footnote""{Both authors were supported in part by grants from the National Science Foundation.}
\footnote""{This paper is in final form and no version of it will be submitted
for publication elsewhere.}

\specialhead 1. Introduction \endspecialhead\medskip

Throughout this paper all rings are commutative, with identity, and
Noetherian, unless otherwise specified.  In particular, ``local ring''
always means Noetherian local ring, unless otherwise specified.
Our objective is to prove a generalization of Faltings' connectedness
theorem [Fal1, Fal2], which asserts that in a complete local domain
$(R,m,K)$ of dimension $n$, if $I\subseteq m$ is an ideal generated by at
most $n-2$ elements, then the punctured spectrum of $R/I$ is connected.  Our
result (see Theorems 3.3 and 3.6) draws the same conclusion without the
hypothesis that $R$ be a domain:  we assume instead that $R$ is complete,
equidimensional (i.e., for every minimal prime $p$ of $R$,
$\dim R/p=\dim R$), and that one of the following conditions, which we shall
prove are equivalent, holds:

a)\quad $H_m^n(R)$ (local cohomology with support in $m$) is indecomposable.

b)\quad The canonical module $\omega$ of $R$ is indecomposable.

c)\quad The $S_2$-ification of $R$ is local.

d)\quad For every ideal $I$ of height two or more, Spec$\ R-V(I)$ is connected.

e)\quad Given any two distinct minimal primes $p$, $q$ of $R$, there is a
sequence of minimal primes $p=p_0,\dots,p_i,\dots,p_r=q$ such that for
$0\le i<r$, the height $(p_i+p_{i+1})\le 1$.

\smallskip
\S2 details the properties of canonical modules for a not necessarily
Cohen-Macaulay ring, as well as the process of $S_2$-ification. By and
large the results of \S2 are known, but in some cases we have not found
a convenient reference.  We note here only that our definition is that
a canonical module $\omega$ for  $(R,m,K)$ is a finitely generated
$R$-module such that $\H(\omega,E)\cong H_m^{\dim R}(R)$, where $E$ is
an injective hull for $K$ over $R$.  The main results are developed in
\S3. \bigskip

\specialhead 2.  Canonical modules and ${\bold S}_{\bold
2}$-ification \endspecialhead \medskip

It will be convenient to have a notation for the ideal that turns out to be
the annihilator of the canonical module.

\demo{{\bf (2.1)} Definition}
If $R$ is a local ring we shall denote by $\J(R)$ the largest ideal which is
a submodule of $R$ of dimension smaller than $\dim R$.  Then $\J (R)$ is nonzero
if and only if some prime $P$ of $\text{Ass } R$ is such that
$\dim R/P<\dim R$, and then $\J (R)\supseteq \text{Ann}_RP$.  Thus,
$\J(R)=(0)$ iff $R$ is equidimensional and unmixed (where {\it unmixed}
means that $(0)$ has no embedded primes).  Moreover, $\J(R)$ consists of
all elements $r\in R$ such that $\dim R/\text{Ann}_R r<\dim R$.
\enddemo

Throughout this section $E=E_R(K)$ denotes an injective hull of the residue
field of the local ring $(R,m,K)$ and ${}^\vee$ denotes the exact functor
$\H(\underline{\phantom M}, E)$ on $R$-modules. We begin by summarizing
many of the known properties of canonical modules, most of which we shall need in
this paper.

\demo{{\bf (2.2)} Remark}
Let $(R,m,K)$ be a local ring with $\dim R=d$.

a)\quad If $R$ is complete, then $R$ has a canonical module, and any
canonical module is isomorphic with $H_m^d(R)^\vee$.

b)\quad Any two canonical modules for $R$ are (non-canonically) isomorphic.

c)\quad If $R$ is a homomorphic image of a Gorenstein ring, then $R$ has a
canonical module.  If $R=S/J$, where $S$ is local, then
$\text{Ext}_S^h(R,S)$ is a canonical module for $R$, where
$h=\dim S-\dim R$.  More generally, if $S\to R$ is local, $R$ is module-finite
over the image of $S$, $S$ is Cohen-Macaulay with canonical module
$\omega_S$, and $h=\dim S -\dim R$, then $\text{Ext}_S^h(R,\omega_S)$ is a
canonical module for $R$.  In particular, if $R$ is a module-finite
extension of a regular (or Gorenstein) local ring $A$, then
$\text{Hom}_A(R,A)$ is a canonical module for $R$.  (The same holds when $R$
is module-finite over the image of $A$ and the two have the same dimension.)

d)\quad A canonical module for $R$ must be killed by $\J(R)$, and is also a
canonical module for $R/\J(R)$, while any canonical module for
$R/\J(R)$ is a canonical module for $R$.  Thus, $R$ has a canonical module
if and only if $R/\J(R)$ has a canonical module.
\medskip
 For parts e)-i) we let $(R,m,K)$ be a local ring with canonical module $\omega$.
\medskip
e)\quad The kernel of the map $R\to\HO$ is $\J(R)$.  Thus, $\omega$ is
faithful if and only if $R$ is equidimensional and unmixed.

f)\quad The module $\omega$ and its completion are both $S_2$.  Moreover,
$\HO$ is a commutative semilocal ring module-finite over the image of $R$
and it is $S_2$ both as an $R$-module and as a ring in its own right.  It
may be identified with a subring of the total quotient ring of $R/\J(R)$.
Moreover, its $m$-adic completion is $S_2$.

g)\quad For every prime $P$ of $R$ such that $\dim R/P=\dim R$, the ring
$(R/P)^\wedge \cong \hat R/P\hat R$ is equidimensional and unmixed.  If
$\J(R)=(0)$ then $\J(\hat R)=(0)$.

h)\quad $R\to\HO$ is an isomorphism if and only if $R$ is $S_2$ and
equidimensional (the latter condition follows from $S_2$ if $R$ is catenary),
and also iff $\hat R$ is $S_2$.  Thus, if $R$ has a canonical module and
$R$ is equidimensional and $S_2$, then $\hat R$ is $S_2$.

i) If $R$ is equidimensional, then for every prime $P$ of $R$, $\omega_P$ is a
canonical module for $R_P$.
\medskip
For parts j)-l) we suppose that $(R,m,K)$ is local with $\J(R)=(0)$ and let $\omega$ be a
canonical module for $R$.  Let $S=\HO$.
Let $m_1,\dots,m_s$ denote
the maximal ideals of $S$.  Note that $\omega$ is an $S$-module, precisely
because $S=\HO$.  Then:
\medskip
j)\quad Every maximal ideal of $S$ has height equal to $\dim R$.

k)\quad When $R$ is complete, so that $S$ is product of local rings $S_i$,
one for every maximal ideal $m_i$ of $S$, and $\omega$ is, correspondingly, a
product of modules $\omega_i$ over the various $S_i$, then $\omega_i$ is
the canonical module for $S_i$ for every $i$.

l)\quad The module $\omega$ is a canonical module for $S$ in the sense
that $\omega_Q$ is a canonical module for $S_Q$ for every prime ideal $Q$
of $S$. \enddemo

\demo{Proof}
Parts a) - c) are standard and can be found in [HK]. The rest of the
results can be found in either [G] or [A1-2]. In particular,
see Theorem 3.2 and Corollary 4.3 in [A2].
\enddemo
\medskip
\demo{{\bf (2.3)}Discussion}
Let $(R,m,K)$ be an equidimensional and unmixed local ring, i.e., such that
$\J(R)=(0)$.  We shall say that a ring $S$ is an $S_2$-ification of $R$ if
it lies between $R$ and its total quotient ring, is module-finite over $R$,
is $S_2$ as an $R$-module, and has the property that for every element
$s\in S-R$, the ideal $\Cal D(s)$, defined as $\{r\in R: rs\in R\}$, has
height at least two.  We are interested in this notion because if
$\J(R)=(0)$ and $R$ has a canonical module $\omega$ then it has an
$S_2$-ification, to wit, $\HO$ identified with a subring of the total quotient
ring of $R$.  Moreover, whenever $R$ has an $S_2$-ification, it is unique.
We prove several propositions in this direction.
\enddemo

\proclaim{(2.4) Proposition}
Let $(R,m,K)$ be a local ring with $\J(R)=(0)$ and let $T$ be its total
quotient ring.  If $f\in T$ let $\Cal D(f)=\{r\in R: rf\in R\}$.  Let
$S$ be the subring of $T$ consisting of all elements $f\in T$ such
that $ht\ \Cal D(f)\ge 2$.  Then $R$ has an $S_2$-ification if and only
if $S$ is module-finite over $R$, in which case $S$ is the
unique $S_2$-ification of $R$.
\endproclaim

\demo{Proof}
It is easy to verify that $S$ is a subring of $T$ containing $R$,
since $\Cal D(r)=R$ for $r\in R$ (and the height is $+\infty$),
$\Cal D(s\pm s')\supseteq \Cal D(s)\cap \Cal D(s')$, and $\Cal D(ss')
\supseteq\Cal D(s)\Cal D(s')$.  Moreover, it is immediate from the way that
we defined an $S_2$-ification that it must be contained in $T$.  We next
observe that if $S_0\subseteq S \subseteq T$ with $S_0$ module-finite over $R$
and $S_0$ is $S_2$ as an $R$-module then $S_0=S$.  To see this, suppose
that $f\in S-S_0$.  Since $\Cal D(f)$ has height at least two (but
cannot be $R$) and $S_0$ is a faithful $R$-module, we must have that there
is a regular sequence $x$, $y$ of length two on $S_0$ in $\Cal D(f)$.  Now
$xf$, $yf\in R$, and so we have that $x(yf)-y(xf)=0$ is a relation on
$x$, $y$ with coefficients in $S_0$.  It follows that $xf\in xS$, so that
$xf=xs$ with $s\in S_0$.  But $x$ is a nonzerodivisor in $S_0$, hence in $R$, and so also in $T$, the total quotient ring of $R$.  Thus, $f=s$, and $f\in S_0$.

Now suppose that $S$ is module-finite over $R$.  We must show that
$S$ is $S_2$.  The depth of $S$ on a height one ideal of $R$ is at
least one, since the ideal must contain a nonzerodivisor of $R$ (we have that
$\J(R)=0$) and this will be a nonzerodivisor in $T$.  Suppose that $I$ has
height at least two.  Choose elements $x$, $y$ in $I$, nonzerodivisors, such
that $(x,y)R$ has height two.  We claim that $x$, $y$ form a regular sequence
on $S$ (and this will complete the proof).  For suppose that we have a
relation $xs=ys'$ with elements $s$ and $s'$ of $S$.  In the total quotient
ring let $f=s/y=s'/x$.  Choose ideals $I$, $I'$ of $R$ of height at least
two such that $Is\subseteq R$ and $I's'\subseteq R$.  Then $II'(x,y)$
multiplies $f$ into $R$ (since $I'xf=I's'$ and $Iyf=Is$), and this ideal has
height two.  Thus, $f\in S$, and so $s'=xf\in xS$. \qed
\enddemo

\demo{{\bf (2.5)} Remarks}
Note that if an $R$-algebra $S'$ is $R$-isomorphic to the $S_2$-ification
$S$ of $R$, then there is a unique $R$-isomorphism $S'\cong S$.  (For each
element $f$ of $S'$ we can choose a nonzerodivisor $r\in R$ such $fr=r'\in R$.
Then if $\phi:S'\to S$ is the isomorphism we must have
$\phi(r')=\phi(fr)=r\phi(f)$, which determines $\phi(f)$ uniquely.)  Thus,
we shall talk about $S_2$-ifications which are not literally subrings of
the total quotient ring of $R$:  they are always, however, canonically
identifiable with such a subring.

Note also that if $S$ is an $S_2$-ification of $R$, then we can choose
finitely many generators for $S$ as an $R$-module, and for each of these
generators an ideal of height at least two in $R$ that multiplies $S$ into
$R$.  It follows that there is an ideal of height at least two in $R$ that
multiplies $S$ into $R$:  intersect the ideals chosen for the individual
generators.
\enddemo

\proclaim{(2.6) Proposition}
If $(R,m,K)$ has an $S_2$-ification $S$ then for every prime ideal $P$ of
$R$, $S_P$ is an $S_2$-ification of $R_P$.
\endproclaim

\demo{Proof}
$S_P$ is $S_2$ over $R_P$, module-finite, and identifiable with a
subring of the total quotient ring.  Any element has a unit multiple of the
form $s/1$, where $s$ is in $S$.  Then $\Cal D(s/1)\subseteq
\Cal D(s)_P\subseteq R_P$ has height at least two. \qed
\enddemo

\proclaim{(2.7) Proposition}
If $(R,m,K)$ is local, $\J(R)=0$, and $\omega$ is a canonical module for
$R$, then $R\to\HO$ is an $S_2$-ification of $R$.
\endproclaim

\demo{Proof}
We know by 2.2 f) that $\HO$ is a module-finite extension of $R$ that may be
identified with a subring of the total quotient ring of $R$.  We also know
that it is $S_2$.  Therefore, it will suffice to show that if
$s\in\HO$, then $\Cal D(s)$ has height at least two.  If not, it will be
contained in some height one prime $P$ of $R$, and $P$ will be in the
support of $R/\Cal D(s)\cong (R(id_\omega)+Rs)/R(id_\omega)\subseteq\HO/R$,
and $\HO/R$ is not supported at any height
one prime by 2.2 h) and i). \qed
\enddemo

\bigskip
\specialhead 3.  Connectedness theorems \endspecialhead \medskip

We first give the statement of the
Faltings' connectedness theorem for complete domains (see [Fal1],
[Fal2]) . In our improvement we will follow the lines of an
argument given in [BR]. We will use the following result, the local
Hartshorne-Lichtenbaum vanishing theorem (see [Ha], [CS], and [BH]).  In its
basic form, it asserts that for a complete local domain $(R,m,K)$ of
dimension $n$, if $I\subseteq m$ is not primary to $m$ then $H_I^n(R)=0$.
In a more precise form, it asserts that if $(R,m,K)$ is any complete local
ring of dimension $n$ and $I\subseteq m$ is an ideal, then
$H_I^n(R)\not= 0$ if and only if there exists a minimal prime ideal $p$
of $R$ such that $\dim R/p=\dim R$ and $I+p$ is primary to $m$.

\proclaim{(3.1) Theorem (Faltings' connectedness theorem)}
Let $(R,m,K)$ be an analytically irreducible local ring of dimension $n$,
and let $\fraktur A$ be an ideal of $R$ generated by at most $n-2$ elements.
Then the punctured spectrum of $\text{Spec } R/\fraktur A$ is connected.
(In other words, there do not exist ideals $I$, $J$ of $R$ such that
$\Rad I\cap J=\Rad \fraktur A$ and $\Rad (I+J)=m$
unless one of the ideals is primary to $m$ and the other has the same radical
as $\fraktur A$.)
\endproclaim

\demo{{\bf (3.2)} Remarks}
We do not even need the condition that $\fraktur A$ be generated by
$n-2$ or fewer elements:  all that is needed is that $H_{\fraktur A}^{n-1}(R)=
H_{\fraktur A}^n(R)=0$, and the second condition is automatic if $R$ is a
complete local domain and $\fraktur A$ is not $m$-primary.

Faltings' original proof was for the equicharacteristic case.
In [BR] a much simpler proof was given, whose outline we shall follow here.
In [HH2], \S6 (see also [HH1]) it is shown that the
integral closure $R^+$ of a complete (or excellent) local domain of positive
characteristic $p$ in an algebraic closure of its fraction field is
``\CM'' in the sense that every system of parameters in $R$ is a regular
sequence in $R^+$.  The characteristic $p$ case of the Faltings'
connectedness theorem can be deduced from the \CM property for $R^+$, and
then the equal characteristic $0$ case also follows by the technique of
reduction to characteristic $p$. It was this point of view which led us to
suspect that the condition that $R$ is a domain could be weakened.
\enddemo

Our next main objective here is to generalize so that $R$ need not be a
domain.  We shall state all of our results in the complete case.  In each
instance, one can achieve the illusion of greater generality by starting
with an arbitrary local ring and requiring that its completion satisfy the
hypotheses we want.

\proclaim{(3.3) Theorem}
Let $(R,m,K)$ be a complete equidimensional local ring such that,
equivalently, $H_m^n(R)$ is an indecomposable $R$-module or such that
$\omega_R$, its canonical module, is indecomposable.  (This is automatic
if $R$ is a domain, since $\omega_R$ is then an ideal of $R$.)  Let
$\fraktur A$ be a proper ideal of $R$ generated by $n-2$ or fewer
elements.  Then the punctured spectrum of $R/\fraktur A$ is connected.
\endproclaim

\demo{Proof}
If not let $I$, $J$ be ideals which give a disconnection, so that $I\cap J$
has the same radical as $\fraktur A$, $I+J$ is primary to $m$, but neither
$I$ nor $J$ is primary to $m$.
The Mayer-Vietoris
sequence for local cohomology then yields:
$$
\cdots\to H_{I\cap J}^{n-1}(R)\to H_{I+J}^n(R)\to H_I^n(R)\oplus
H_J^n(R)\to H_{I\cap J}^n(R)\to\cdots
$$
and the first and last terms displayed are zero, since $I\cap J$ has the same
radical as an ideal with at most $n-2$ generators.  Since $I+J$ is primary
to $m$, this yields an isomorphism:
$$
H_m^n(R)\cong H_I^n(R)\oplus H_J^n(R).
$$
The fact that $H_m^n(R)$ is indecomposable implies that one of
the summands, say $H_J^n(R)$, is zero.  But then the local
Hartshorne-Lichtenbaum vanishing theorem implies that for every prime
$P\in\Cal P$, the set of minimal primes of $R$, $J+P$ is not primary to $m$.
Let $P$ be one of these minimal primes. The intersection of $I+P$ and
$J+P$ is still, up to radicals, $\fraktur A+P$, while the sum is still
primary to $m$.  Thus, applying the local Hartshorne-Lichtenbaum
vanishing theorem to the domain $R/P$, we see that $I+P$ must be primary
to $m$ for every minimal prime $P$.  But this implies that
$\underset{P\in\Cal P}\to\bigcap (I+P)$ is primary to $m$, and up to radicals this is the
same as $I+ \underset{P\in\Cal P}\to\bigcap P$.  Since $\underset{P\in\Cal P}\to\bigcap P$
is the ideal of nilpotents, we find that $I$ itself is primary to $m$, a
contradiction. \qed
\enddemo

This result motivates a study of when the canonical module of a complete
local equidimensional ring is indecomposable.  We begin by associating a
graph with such a ring.

\demo{{\bf (3.4)} Definition}
Let $R$ be an equidimensional local ring.  We denote by $\Gamma_R$ the
(undirected) graph whose vertices are the minimal primes of $R$, and whose
edges are determined by the following rule:  if $P$, $Q$ are distinct
minimal primes of $R$, then $\{P,Q\}$ is an edge of $\Gamma_R$ if and only
if $P+Q$ has height one.
\enddemo

We next observe:

\proclaim{(3.5) Proposition}
Let $(R,m,K)$ be a local ring with canonical module $\omega$.  Suppose that
$\J(R)=(0)$, and let $S=\HO$, the $S_2$-ification of $R$.  Then:

a)\quad For every prime ideal $P$ of $R$ and $Q$ of $S$, if $Q$ lies over
$P$ then $ht\ Q=ht\ P$.

b)\quad For every ideal $I$ of $R$, $\text{height }IS=\text{height } I$.

c)\quad Contraction gives a bijection of the minimal primes of $S$ with
the minimal primes of $R$ and a bijection of the height one primes of $S$
with the height one primes of $R$.
\endproclaim

\demo{Proof}
a) To study the primes of $S$ lying over $P$, we first replace $R$, $\omega$,
and $S$, by $R_P$, $\omega_{P}$, and $S_P$.  Thus, there is no loss of
generality in supposing that $P$ is the maximal ideal of $R$.  But then the
result follows from 2.2 j).

b)\quad If $P$ is a prime ideal containing $I$ whose height is the same
as that of $I$, then there is a prime ideal $Q$ of $S$ lying over $P$.
Then $Q$ contains $IS$, and so $ht\ IS\le ht\ Q = \ ht\ P= ht\ I$.  On the other
hand, if $Q$ is a prime ideal containing $IS$ whose height is the same as
that of $IS$, then $Q$ contracts to a prime ideal $P$ containing $I$, and
so $ht\ IS=ht\ Q= ht\ P\ge ht\ I$.

c)\quad Any height $k$ prime $P$ of $R$ has at least one prime of $S$ lying
over it, and all such primes have height $k$ by part a).  Moreover, all
height $k$ primes of $S$ lie over height $k$ primes of $R$.  Thus, it
suffices to show that when $k=0,1$, there is at most one prime of $S$ lying
over $P$.   But the primes of $S$ lying over $P$ correspond to the primes
of $S_P$ lying over $PS_{P}$, and $S_P$ is the $S_2$-ification of $R_P$.
When $\dim R_P\le 1$, $R_P$ is its own $S_2$-ification, and the result
follows. \qed
\enddemo

We are now ready to prove a central result:

\proclaim{(3.6) Theorem}
Let $(R,m,K)$ be a complete local equidimensional ring with $\dim\ R=n$.  The the following conditions are equivalent:

a)\quad $H_m^n(R)$ is indecomposable.

b)\quad The canonical module $\omega=\omega_R$ of $R$ is indecomposable.

c)\quad The $S_2$-ification $S$ of $R/\J(R)$ is local.

d)\quad For every ideal $J$ of height at least two, $\text{Spec } R-V(J)$
is connected.

e)\quad $\Gamma_R$ is connected.
\endproclaim

\demo{Proof}
We shall prove that $a)\Leftrightarrow b) \Leftrightarrow c)\Rightarrow
d)\Leftrightarrow e)\Rightarrow c)$.  The equivalence of a) and b) is
clear.  Now assume b).  The module $\omega$ is also a canonical module
for $R/\J(R)$ and, consequently, for the $S_2$-ification
$\HO=\HH_{R/\J(R)}(\omega, \omega)$ of $R/\J(R)$ as well.  If the $S$
is not local, then $\omega$ is a product of nonzero factors
corresponding to the various factors rings of $S$, and this will yield
a non-trivial direct sum decomposition of $\omega$ over $R$.  On the
other hand, if $S\cong\HO$ is local, it contains no idempotents other
than $0$, $1$, and this implies that $\omega$ is indecomposable.  Thus,
$c)\Rightarrow b)$ as well, and we have established the equivalence of
the first three conditions.

Now assume that $S$ is local, and that $I$, $I'\subseteq m$ are such
that $I\cap I'$ is nilpotent but $I+I'$ has height at least two.  We
can replace $I$, $I'$ by powers and assume that $II'=0$ but $I+I'$ has
height at least two.  This situation is preserved when we pass to
$R/\J(R)$, for when $R$ is equidimensional $\J(R)$ consists of
nilpotents.  Thus, we might as well assume that $\J(R)=(0)$.  Next,
note that $IS+I'S$ has height at least two (since its height is the
same as that of $I+I'$, by the preceding proposition).  Moreover, if
neither $I$ nor $I'$ is primary to $m$ (i.e., if neither has height
equal to $\dim R$) then neither $IS$ nor $I'S$ is primary to the
maximal ideal of $S$.

Thus, there is no loss of generality in assuming that $R$ is $S_2$. The ideal
$I+I'$ will contain a regular sequence $u+u'$, $v+v'$ of length two, where
$u$, $v\in I$ and $u'$, $v'\in I'$.  The relation $v(u+u')-u(v+v')=0$ then
shows that $u\in (u+u')R$, while $u'\in (u+u')R$ similarly.  This yields
$u=(u+u')a$, $u'=(u+u')b$, and so $u+u'=(u+u')(a+b)$.  Since $u+u'$ is not
a zerodivisor, $1=a+b$, and it follows that at least one of $a,b$ is a unit.
Suppose that $a$ is a unit:  the other case is similar.  Then
$u=(u+u')a$ implies that $u$ is a nonzerodivisor, while $uI'\subseteq II'=(0)$
then implies that $I'=0$, and so $I$ is primary to $m$.  This completes the
proof that $c)\Rightarrow d)$.

We next want to see that $d)\Leftrightarrow e)$.  Suppose that one has
ideals $I$, $I'$ such that $I\cap I'$ is nilpotent. Then we can replace
$I$, $I'$ by their radicals while only increasing $I+I'$.  Then each of
$I$, $I'$ is a finite intersection of primes.  For each minimal prime
$p$ of $R$, $p\supseteq I\cap I'$, and so $p$ must contain either a
minimal prime of $I$ or a minimal prime of $I'$.  Thus, $p$ must be
either a minimal prime of $I$ or a minimal prime of $I'$.  If we omit
all non-minimal primes from the primary decomposition of $I$
(respectively, $I'$) and intersect the others, we get two larger ideals
whose intersection is still $\Rad (0)$.  Thus, it is possible to give
$I$, $I'$ such that $\Rad (I\cap I')=\Rad (0)$ and $I+I'$ has height
two if and only if one can do this with ideals $I$, $I'$ coming from a
partition of the minimal primes of $R$ into two nonempty sets, with $I$
the intersection of the minimal primes in one set and $I'$ the
intesection of the minimal primes in the other set.  If one set
consists of $\{p_1,\dots,p_h\}$ and the other of $\{q_1,\dots,q_k\}$ we
shall have $I=\cap_i p_i$, $I'=\cap_j q_j$, and $I+I'$ will then have
the same radical as $\cap_{i,j}(p_i+q_j)$, and will have height at
least two if and only if every $p_i+q_j$ has height at least two.
Thus, d) fails if and only if the minimal primes can be partitioned
into two nonempty sets such that no edge of $\Gamma$ joins a vertex in
one set to a vertex in the other,  which is precisely the condition for
$\Gamma$ to be disconnected.  Thus, $d\Leftrightarrow e)$.

Finally, we show that $e)\Rightarrow c)$.  Suppose that $\Gamma$ is
connected.  We want to prove that the $S_2$-ification of $R/\J(R)$ is
local.  The graph associated with $R/\J(R)$ is the same as that
associated with $R$, so that we may assume that $\J(R)=(0)$.  If the
$S_2$-ification $S$ of $R$ has two or more maximal ideals, say $\Cal
M_1,\dots,\Cal M_r$, where $r\ge 2$, for each $\Cal M_j$ let $\Cal P_j$
denote the set of minimal primes of $S$ contained in $\Cal M_j$.  Then
$\Cal P_j$ is evidently non-empty.  There is a bijection between the
minimal primes of $S$ and those of $R$, so that for each $\Cal P_j$
there is a corresponding set of minimal primes $\Cal Q_j$ of $R$.  To
complete the argument, it will suffice to show that if $i$, $j$ are
different then it is impossible to have an edge joining a vertex in
$\Cal Q_i$ to a vertex in $\Cal Q_j$.  If there were such an edge,
there would be a height one prime $P$ of $R$ containing both a minimal
prime in $\Cal Q_i$ and a minimal prime in  $\Cal Q_j$.  Then $R_P\cong
S_P$, and it follows that the unique prime of $S$ lying over $P$
contains both a prime of $\Cal P_i$ and a prime of $\Cal P_j$.  Let
$\Cal M$ be a maximal ideal of $R$ containing $P$.  Then $\Cal M$
contains both a prime of $\Cal P_i$ and a prime of $\Cal P_j$, which is
impossible:  $S$ is a finite product of local rings, and each prime
ideal of $S$ is therefore contained in a unique maximal ideal of $S$,
forcing $\Cal M_i=\Cal M=\Cal M_j$. \qed \enddemo

When $R$ is not complete, it is necessary to study the graph associated with
the minimal primes in the completion:  the domain property is frequently
lost when one completes.  However, the characterization in (3.6c) behaves
better, as we see in (3.7) below.  First note that in the sequel, if $(R,m,K)$ is equidimensional (but possibly has $\J(R)\not= (0)$), by an $S_2$-ification
of $R$ we mean an $S_2$-ification for $R/\J(R)$.

\proclaim{(3.7) Corollary}
Let $(R,m,K)$ be an equidimensional local ring with canonical module
$\omega$.  Let $n=\dim R$.  Then $H_m^n(R)$ is indecomposable if and only
if the $S_2$-ification of $R$ is local.  In particular, if $R$ is $S_2$,
then $H_m^n(R)$ is indecomposable.
\endproclaim

\demo{Proof}
Killing $\J(R)$, if it is not zero, does not affect either issue, and so we may
assume that $\J(R) = (0)$.  The $S_2$-ification of $\HR$ is $\HH_{\HR}(\Ho,\Ho)\cong\HR\otimes_R\HO$. Equidimensionality is preserved by completion here, and the issue of whether a semilocal ring is local is not affected by completing with respect to its Jacobson radical. \qed
\enddemo

The next Proposition is well-known (e.g. see [B,5.2] which gives a much
more general result), but we include a proof as it is fairly short.

\proclaim{(3.8) Proposition}
If $(R,m,K)$ is an excellent local, equidimensional ring then $R$ has an
$S_2$-ification $S$, and $\HS$ is the $S_2$-ification of $\HR$.
\endproclaim

\demo{Proof}
We may first kill $\J(R)$, and we henceforth suppose that it is $0$.
Because $R$ is excellent, this is preserved by completion.  Let $S$ be
the set of elements of the total quotient ring $\Cal T(R)$ of $R$ that
are multiplied into $R$ by an ideal of height two or more.  It will
suffice to show that $S$ is module-finite over $R$.  Note that $\Cal
T(R)\subseteq \Cal T(\HR)$, since $\HR$ is flat over $R$.  If $S$ is
not module-finite over $R$ we can choose a sequence of elements
$\{s_i\}$ in $S$ such that the sequence of $R$-submodules $\sum_{i=1}^j
Rs_i\subseteq S$ is strictly increasing with $j$.  Clearly, each
$s_i\in\Cal T(\HR)$.  Since $\HR$ has an $S_2$-ification, we can choose
$j$ so large that $\sum_{i=1}^j \HR s_i =\sum_{i=1}^{j+1} \HR s_i$. We
can choose a nonzerodivisor $a\in R$ such that for $0 \leq i \leq
j+1$,  $as_i \in R$. Then $\sum_{i=1}^j \HR as_i =\sum_{i=1}^{j+1} \HR
as_i$.  Since $\HR$ is faithfully flat over $R$, we find that
 $\sum_{i=1}^j R as_i =\sum_{i=1}^{j+1} R as_i$, and so $as_{j+1}\in
\sum_{i=1}^jRas_i$. Since $a$ is a nonzerodivisor in $R$ (and hence in
$S$), it follows that $s_{j+1}\in\sum_{i=1}^jRs_i$, a contradiction.

Thus, $R$ has an $S_2$-ification, $S$.  When we complete, since the fibers
of $R\to \HR$ are \CM, we see that $\HS$ is $S_2$ as an $\HR$-module.  It is
clear that it is contained in the total quotient ring of $\HR$: $S/R$ is
killed by a nonzerodivisor $a$ in $R$, and so $\HS/\HR$ is also killed by $a$.
Moreover, if $I$ is an ideal of $R$ of height at least two killing $S/R$,
then $I\HR$ kills $\HS/\HR$.  It follows that $\HS$ is the $S_2$-ification
of $\HR$. \qed
\enddemo

\proclaim{(3.9) Proposition}
Let $(R,m,K)$ be an excellent equidimensional local ring.

a)\quad The $S_2$-ification of $R$ is local if and only if the $S_2$-ification
of $\HR$ is local.

b)\quad The $S_2$-ification of $R$ is local if and only if the $S_2$-ification
of $R_{red}$ is local.

c)\quad If $R$ is $S_2$ and $x_1,\dots,x_k$ is a part of a system of
parameters, then the $S_2$-ification of $R/(x_1,\dots,x_k)R$ is local.
\endproclaim

\demo{Proof}
a) This is immediate from (3.8).

b)\quad Since $R$ is excellent, $(\HR)_{red}\cong (R_{red})^\wedge$. Thus,
we may assume that $R$ is complete.  We may also assume that
$\J(R)=0$.  The result then follows from the fact that the graphs
associated with $R$ and $R_{red}$ are the same.

c)\quad The issues are unaffected by completing $R$ and killing the
nilpotents.  It is easy to see that $B=R/(x_1,\dots,x_k)R$ is again
equidimensional.  If the $S_2$-ification is not local then there is a
localization of $B$ at a prime of height at least two such that the punctured
spectrum is disconnected, and this ring may be viewed as a quotient of a
localization $R_Q$ of $R$.  But $R_Q$ is $S_2$ and has a canonical module
(it is a localization of a complete ring, and so is a homomorphic image of a
localization of a complete regular local ring), so that Corollary 3.7 applies.
The result now follows from Theorem 3.3. \qed
\enddemo

\demo{{\bf (3.10)} Remark}
Even for complete domains, the fact that the $S_2$-ification of the local
ring $R$ is local does not imply, in general, that the $S_2$-ification of
every local ring of $R$ is local.  Some primes of $R$ may have more than one
prime of the $S_2$-ification lying over them.

For example, consider $R=K[x,y,yz,z(z-x),z^2(z-x)]\subseteq S=K[x,y,z]$.
This extension is integral, since $z$ satisfies $Z^2-xZ-z(z-x)=0$.  The
element $z$ is multiplied into $R$ by the height two ideal $(y,z(z-x))$.
Now $P=(y,yz,z(z-x),z^2(z-x))R\subseteq R$ is a height two prime, but two
prime ideals of $S$ lie over it:  $(y,z)S$ and $(y,z-x)S$.

If we complete $R$, $S$ at their homogeneous maximal ideals both rings remain
domains.  $\HS$ is the $S_2$-ification of $\HR$, and is local.  However,
there are two primes of $\HS$ lying over $\hat P=P\HR$, and so the
$S_2$-ification of $\HR_{\hat P}$ is not local.
\enddemo

\demo{{\bf (3.11)} Remark}
Suppose that one is trying to give an elementary proof of Faltings'
connectedness theorem in the generality we have obtained here, perhaps
without using local cohomology.  It would suffice to prove that if $R$
is a complete reduced $S_2$ local ring and $x$ is a single parameter, then
the $S_2$-ification  of $R/xR$ is local.  The connectedness theorem can
be reduced to the case of parameters, and if one knows the single fact
stated above, one can carry through an induction on the number of parameters.
However, it is quite possible that the case of a single parameter is no
easier than the general case.  We next note:
\enddemo

\proclaim{(3.12) Proposition}
If $(R,m,K)$, $(S,n,K)$ are two complete equidimensional local rings with
algebraically closed coefficient field $K$ and the $S_2$-ifications of $R$,
$S$ are local then so is the $S_2$-ification of $T=R\hat\otimes_KS$.
\endproclaim

\demo{Proof}
We may assume that $R$, $S$ are reduced, and so is $T$.  We also note that
when $R$, $S$ are domains, then $T$ is a domain.  It follows that every
minimal prime of $T$ has the form $p\hat\otimes_K S+ R\hat\otimes_K q$
where $p$ is a minimal prime of $R$ and $q$ is a minimal prime of $S$.
(Any prime of $T$ will contract to some prime $P$ of $R$, and also to some
prime $Q$ of $S$.  Hence, it must contain $P\hat\otimes_K S+ R\hat\otimes_K Q$,
which is the kernel of the map $R\hat\otimes_K S\twoheadrightarrow
R/P\hat\otimes_K S/Q$.  This immediately shows that the minimal primes are a
subset of the ideals $p\hat\otimes_KS+R\hat\otimes_Kq$ for $p$ minimal in
$R$ and $q$ minimal in $S$.  Since it is easy to see that these ideals are
mutually incomparable, they are all minimal primes.)  Let $\Gamma$, $\Gamma'$
be the graphs associated with $R$, $S$, respectively.  The vertices of the
graph associated with $R\hat\otimes_KS$ are in bijective correspondence with
the set $\Gamma\times\Gamma'$.  There is an edge from $(p,q)$ to
$(p',q')$ if and only if $\text{ht}\left((p+p')\hat\otimes_KS+
R\hat\otimes_K(q+q')S\right)=1$, which happens iff either
$\text{ht}(p+p')=1$ and $q=q'$ or $p=p'$ and $\text{ht}(q+q')=1$.  But then,
in the graph associated with $R\hat\otimes_KS$, with its vertices identified
with $\Gamma\times\Gamma'$, we have that each subgraph $\Gamma\times\{q\}$ is
connected for every $q\in\Gamma'$, and each subgraph $\{p\}\times\Gamma'$
is connected for every $p\in\Gamma$.  It follows that $\Gamma\times\Gamma'$
is connected. \qed
\enddemo

\demo{{\bf (3.13)} Remarks on the graded case}
Now suppose that $R$ is a finitely generated ${\Bbb N}$-graded $K$-algebra
with $R_0=K$.  Let $m$ be the homogeneous maximal ideal of $R$.  If $R$ is
a domain, so is its completion (with respect to the homogeneous maximal
ideal), since its completion has a filtration with respect to which the
associated graded ring is $R$, which is a domain.  This implies that the
homogeneous primes of $R$ remain prime when we complete.  The minimal primes
(in fact, all associated primes) of $R$ are homogeneous, and so correspond
to the minimal primes of the completion.  If $R$ is equidimensional one can
check whether $H_m^n(R)$ is indecomposable by checking whether the graph
associated with the minimal primes of $R$ is connected:  It is not necessary
to complete, since the completion will have the same graph.  One then gets
an immediate family of corollaries of the connectedness theorems given here
for intersections of projective varieties.  One can also apply the technique
of reduction to the diagonal to prove theorems:  it may be desirable in that
case to assume that the field is algebraically closed, so that products of
irreducible components remain irreducible.

It is also worth noting that one can give a graded resolution of $R$ over
a polynomial ring.  Using $\text{Ext}$ to compute the canonical module then
produces a graded canonical module.  The automorphisms of it yield a
``global'' $S_2$-ification $S$ of $R$.  However, while the ring $S$ is a
graded module over $R$, it need not have the property that $S_0$ is $K$:
when $S$ decomposes, one has non-trivial idempotents in $S_0$. \enddemo

\enddocument